\documentclass[a4paper,twoside,10pt]{amsart}
\usepackage{graphicx}

\newtheorem{definition}{Definition}[section]

\newtheorem{remark}{Remark}[section]
\newtheorem{remarks}{Remarks}[section]
\newtheorem{ex}{Example}[section]

\newtheorem{theorem}{Theorem}[section]

\newtheorem{proposition}{Proposition}[section]

\newlength{\szer}
\newcommand{\Teiss}[2]{%
\settowidth{\szer}{$\displaystyle\frac{#1}{#2}$}%
\setlength{\szer}{0.5\szer}%
\left\{\hspace{\szer}%
\raisebox{0.14ex}{\makebox[0pt]{$\displaystyle\frac{#1}{\phantom{#2}}$}}%
\raisebox{-0.14ex}{\makebox[0pt]{$\displaystyle\frac{\phantom{#1}}{#2}$}}%
\hspace{\szer}\right\}%
}
\def\B{{\mathbf B}}
\def\P{{\mathbf P}}
\def\Pp{{\mathcal P}}

\def\R{{\mathbf R}}

\def\C{{\mathbf C}}
\def\Ii{{\mathcal I}}

\def\Newton{{\mathcal N}}

\def\Z{{\mathbf Z}}
\def\O{{\mathcal O}}

\begin{document}
\title{Some resonances of \L ojasiewicz inequalities}
\author{Bernard Teissier}
\address{
Institut Math\'{e}matique de Jussieu, UMR 7586 du CNRS, 175 Rue du Chevaleret, 75013 Paris, France}

\email{teissier@math.jussieu.fr}

\maketitle
\hfill{\it To the memory of Stanis\l aw \L ojasiewicz}\par\medskip\noindent
\vskip.2truecm
\centerline{Summary\footnote{A.M.S 1991 Classification: 32S15 (14B05 14P15)}}\par\medskip\noindent
\L ojasiewicz's inequalities were born to solve a problem in analysis (see \cite{L} and \cite {Ma}) and, as deep mathematical ideas will, they resonate in other fields. I chose among a very large number of possibilities to try to survey three such resonances in commutative algebra and algebraic geometry. The first section summarizes the interpretation of the best possible exponent $\theta$ for a local \L ojasiewicz inequality of an holomorphic function $g$ with respect to an ideal generated by holomorphic functions $(f_k)_{1\leq k\leq r}$, written as $\vert g(z)\vert^\theta\leq C{\rm sup}_k\vert f_k(z)\vert$, as the inclination of an edge of a Newton polygon associated to the "dicritical" components of a log resolution of the ideal. The second calls attention to recent results which show that some rational numbers cannot be \L ojasiewicz exponents for the gradient inequality in the plane, and the final one reports on a recent result of Moret-Bailly which opens perspectives for a \L ojasiewicz inequality in infinite dimensional spaces. No new result is presented here and I have left out a considerable body of significant results directly related to the computation, interpretation or generalization of the idea of a \L ojasiewicz inequality.\par
\emph{Stanis\l aw \L ojasiewicz was a remarkable mathematician in many ways, and it is a privilege to have shared a friendship with him from our meeting in Carg\`ese in 1972 to his death and to be able to honor his memory.}
\section{\L ojasiewicz exponents and Newton polygons}
We consider the set of local Newton polygons in the first quadrant of $\R^2$. It is the set of boundaries of the convex hulls of finite unions of subsets of the form $(a,b)+\R_+^2$ where $(a,b)\in \R_+^2$. There is a natural commutative operation, the sum, on Newton polygons. It is obtained by taking the boundary of the Minkowski sum of the infinite convex regions which two polygons bound. It is commutative and associative. Its unit element is the Newton polygon consisting of the union of the two coordinate half-lines.\par
The following notation, introduced in \cite{T1}, will be convenient for us: given two positive numbers, which for us will be integers, denote by $\Teiss{\ell}{h}$ the elementary Newton polygon which is the boundary of the convex hull of the subset $((0,h)+\R_+^2)\bigcup ((\ell, 0)+\R_+^2)$ of $\R^2$.\par
If we add to these finite elementary Newton polygons the two elementary polygons corresponding to the case where $h$ or $\ell$ (but not both) may be taken to be infinite, we obtain a set of generators for the monoid of all Newton polygons. Every Newton polygon may be written as 
$$\Newton= \sum_{i=1}^s\Teiss{\ell_i}{h_i},$$
and if we require that the \emph{inclinations} $\frac{\ell_i}{h_i}$ be all different, this decomposition is unique. The unit element can be written as $\Teiss{0}{0}$. See \cite{T1} or \cite{L-T}, Compl\'ement 2, for details.\par\noindent

\vspace{1in}
\begin{picture}(2,1)\setlength{\unitlength}{.4cm}
\put(0,0){\line(0,1){5}}
\put(0,0){\line(1,0){5}}
\put(0,3){\line(4,-3){4}}
\put(-0.9,3){$h$}
\put(4,-0.7){$\ell$}
\put(6, 3){$=$}
\put(7,3){$\Teiss{\ell}{h}$}
\put(14,0){\line(0,1){5}}
\put(14,0){\line(1,0){6}}
\put(18,0){\line(0,1){5}}
%\put(13,6){$h_2$}
\put(18,-0.7){$\ell'$}
\put(22, 3){$=$}
\put(23,3){$\Teiss{\ell'}{\infty}$}

\put(0,-6){$\Teiss{\ell}{h}+\Teiss{\ell'}{\infty}$}
\put(7,-6){$=$}
\put(9,-7){$h$}
\put(10,-10){\line(0,1){8}}
\put(10,-10){\line(1,0){10}}
\put(14, -7){\line(0,1){5}}
\put(14,-7){\line(4,-3){4}}
\put(18,-11){$\ell+\ell'$}
\put(14,-11){$\ell'$}
\end{picture}
\vspace{2in}

\par Let now $(X,\O_X)$ be a reduced equidimensional complex analytic space, and $\Ii$ a sheaf of ideals defining a nowhere dense analytic subspace $Y\subset X$. The normalized blowing up $\pi\colon NB_{\Ii}X\to X$ of $\Ii$ is a proper and bimeromorphic analytic map, where the space $NB_{\Ii}X$ is normal and the ideal $\Ii\O_{NB_{\Ii}X}$ is invertible.\par
Given a point $x\in X$ let us consider the finitely many irreducible components $(D_i)_{1\leq i\leq t}$ of the exceptional divisor of $\pi$ which meet the compact fiber $\pi^{-1}(x)$. For each component we may consider the general fiber of the map $\pi\colon \vert D_i\vert\to \vert\pi( D_i)\vert$ and the underlying reduced general fiber. Its degree with respect to the $\pi$-ample sheaf $\Ii\O_{NB_{\Ii}X}/\Ii^2\O_{NB_{\Ii}X}$ will be denoted by ${\rm deg}_\pi \vert D_i\vert$.\par\noindent
Now to any germ of analytic function $g\in \O_{X,x}$ we can associate the \emph{order $v_{D_i}(g\circ \pi)$ of vanishing} \emph{along $\vert D_i\vert$} of the composition with $\pi$ of a representative of $g$ at a general point of $\vert D_i\vert$. We may replace $g$ by an ideal $J$ if we define the order of vanishing $v_{D_i}(J\circ \pi)$ as the infimum of the orders of vanishing of the (compositions with $\pi$ of) elements of $J$, which is attained by some generators of the ideal.\par\noindent
We can now define the \emph{Newton polygon of $g$ with respect to $\Ii$ at $x$} as
$$\Newton_{\Ii,x}(g)=\sum_{i=1}^t{\rm deg}_\pi \vert D_i\vert\Teiss{v_{D_i}(\Ii)}{v_{D_i}(g)},$$
provided that we have $g\in\sqrt{\Ii_x}$, which ensures that none of the $v_{D_i}(g)$ is zero, and the \emph{Newton polygon of $J$ with respect to $\Ii$} as
$$\Newton_{\Ii,x}(J)=\sum_{i=1}^t{\rm deg}_\pi \vert D_i\vert\Teiss{v_{D_i}(\Ii)}{v_{D_i}(J)},$$
provided that we have $J\subset\sqrt{\Ii_x}$.
\begin{remarks} 1) Each element $g\in\sqrt{\Ii_x}$ determines a pre-ordering on the set of $D_i$, by reading from left to right the compact faces of $\Newton_{\Ii,x}(g)$. It would be interesting to determine, for a given $\Ii$ and $x$, which pre-orderings are realizable in this manner. For example, if the germ $\Ii_x$ is a primary ideal for the maximal ideal $m_x$, the trivial pre-ordering, where all elements are equal, is realizable by a superficial element $g$ for $\Ii_x$ \emph{(see \cite{L-T}, Compl\'ement 2)}.\par\noindent
2) The strict transforms of the components $D_i$ in a log-resolution of the pair $(X,Y)$ may be called the \emph{dicritical} components at $x$ of the exceptional divisor of the resolution. In the special case where $\Ii$ is a primary ideal for the maximal ideal of a point on a normal surface $X$, the normalized blowing up of $\Ii$ is equal to the normalized blowing-up of the ideal generated by two general elements $(f,g)$ of $\Ii$ and the strict transforms of the $D_i$ are exactly the components of the resolution where the quotient $\frac{\tilde f}{\tilde g}$ of the strict transforms of $f$ and $g$ is not constant. The computation of normalized blowing-ups is difficult in dimension $>2$, although there exist algorithms to compute the normalization.\par\noindent A method of computation of the \L ojasewicz exponent in the case of an ideal primary for the maximal ideal of $\O_{X,x}$, through an explicit log resolution is explained in \cite{BA-E}.\par\noindent
3) In the case where the ideal $\Ii_x$ is primary for the maximal ideal of $\O_{X,x}$, so that all the components $D_i$ map to $x$, the valuations $v_{D_i}$ are called the \emph{Rees valuations} associated to that ideal \emph{(see \cite{R1})}.
\end{remarks}
What we call "the" \L ojasiewicz exponent is of course the smallest real number for which the inequality holds in some neighborhood of a given point.\par
\begin{proposition}{\rm (see \cite{L-T}, Th\'eor\`eme 4.6).} Given a set of generators $f_1,\ldots, f_s$ for the germ $\Ii_x\in \O_{X,x}$ and $g\in \O_{X,x}$, the smallest real number $\theta$ such that there exist a neighborhood $U$ of $x$ in $X$ and a constant $C>0$ such that the inequality
$$\vert g(z)\vert^\theta\leq C{\rm sup}_{k=1}^s \vert f_k(z)\vert $$
holds for all $z\in U$, is equal to $$\theta={\rm max}_{i=1}^t \frac{v_{D_i}(\Ii)}{v_{D_i}(g)}.$$
It is therefore a rational number. If $g\notin\sqrt{\Ii_x}$, set $\theta=+\infty$.
\end{proposition}
For the rationality of the \L ojasiewivz exponent in the real-analytic case see Risler's appendix in \cite{L-T}, \cite{B-R}, and \cite{Fe}.\par\medskip\noindent
Given $g\in\sqrt{\Ii_x}$, denote by $\nu_{\Ii_x}(g)$ the largest integer $n$ such that $g\in \Ii_x^n$.
In \cite{L-T}, we prove that the limit $\overline\nu_{\Ii_x}(g)={\rm lim}_{k\to\infty}\frac{\nu_{\Ii_x}(g^k)}{k}$ exists and is equal to ${\rm min}_{i=1}^t\frac{v_{D_i}(g)}{v_{D_i}(\Ii)}$. Thus, the best \L ojasiewicz exponent, which is the largest inclination of the compact edges of $\Newton_{\Ii,x}(g)$, can be understood as the inverse of this limit. We also give some indications on the history of these results.\par
Again, one can extend the definition of $\overline \nu$ to an ideal, by defining $\nu_{\Ii_x}(J)$ to be the largest integer such that $J\subset \Ii_x^n$ and proceeding as above.\par
The fact that the \L ojasiewicz exponents appears as (the inverse of) the inclination of a compact edge of a Newton polygon\footnote{There are at least two other such constructions for $\overline\nu_{\Ii}(g)$: one appears in \cite{L-T}, appendice au \S 4, and the other in \cite{Hi3}.}  is interesting in particular because other features of the same polygon have, at least in special cases, algebro-geometric interpretations:\par\noindent
$\bullet$ If the ideal $\Ii_x$ is primary for the maximal ideal, the length of the horizontal projection of $\Newton_{\Ii_x}(g)$ is equal to the Samuel multiplicity $e(\Ii_x)$ of the primary ideal $\Ii_x$ in $\O_{X,x}$.\par\noindent
$\bullet$ The length of the vertical projection is equal to the \emph{degree function} of $g$ with respect to $\Ii_x$, studied by David Rees and which is the Samuel multiplicity of the image of $\Ii_x$ in the quotient $\O_{X,x}/g\O_{X,x}$. In particular if $g$ is a general element of another primary ideal $J$, it is the mixed multiplicity $e(\Ii_x^{[d-1]},J^{[1]})$.\par\noindent
$\bullet$ As an immediate consequence of the convexity of the Newton polygon, setting $d={\rm dim}X$, we have the inequality $\overline\nu_{\Ii_x}(J)\leq \frac{e(\Ii_x^{[d-1]},J^{[1]})}{e(\Ii_x)}$.\par\noindent

\begin{subsection}{\L ojasiewicz exponent and plane sections}
In the paper \cite{Hi3}, Hickel has shown that given a primary ideal $\Ii$ in a regular local ring $(\O,m)$ of dimension $d$ and equal characteristic zero, if we study the $\overline\nu_\Ii(m)$, the $\overline\nu_{\Ii\O_{H^{(i)}}}(m\O_{H^{(i)}})$ is well defined for a quotient $\O_{H^{(i)}}$ of $\O$ by the ideal generated by $d-i$ "general" linear forms. This suggests the following:\par\medskip\noindent
\textbf{Problem}: Let $\O$ be a Cohen-Macaulay excellent equicharacteristic normal local ring of dimension $d$ and let $\Ii$ and $J$ be two primary ideals. Given a family ${\mathbf g}=(g_1,\ldots ,g_{d-i})$ of $d-i$ elements of $J$, let us denote by $J_{\mathbf g}$ the ideal of $\O$ which they generate. Is it true that the Newton polygon $\Newton_{\Ii\frac{\O}{J_{\mathbf g}}}(J\frac{\O}{J_{\mathbf g}})$ is independent of the family ${\mathbf g}$ provided it is "sufficiently general"?\par\medskip
 
Such a result is true at least in the case where $\Ii\subset \C\{z_1,\ldots ,z_d\}$ is the jacobian ideal of a hypersurface $f(z_1,\ldots ,z_d)=0$ with an isolated singularity at the origin i.e., the ideal generated by the partial derivatives, and $J$ is the maximal ideal $m$. The Newton polygon $\Newton_{j(f)}(m)$ then takes the name of \emph{jacobian Newton polygon}. There are four ingredients:\par\noindent
$\bullet$ One proves that the jacobian Newton polygon of the hypersurfaces of an analytic Whitney equisingular family is constant (see \cite{T2}). \par\noindent
$\bullet$ The general hyperplane sections of an analytic space form a Whitney equisingular family (see \cite{T2}).\par\noindent
$\bullet$ The image in the local ring of a hyperplane of the jacobian ideal of $f$ and the jacobian ideal of the restriction of $f$ to the hyperplane have the same integral closure if the hyperplane is "general"  (a special case of the \emph{idealistic Bertini Theorem}; see \cite{T1}, 2.8).\par\noindent
$\bullet$ The Newton polygon depends only on the integral closures of the ideals (see \cite{L-T}). 
Here a hyperplane is "general" if it is not a limit position of tangent hyperplanes to the hypersurface along sequences of nonsingular points tending to the singular point.\par
The way to extend this argument to the general case is not clear.
\end{subsection}

\begin{subsection}{The \L ojasiewicz exponent at infinity}
The \L ojasiewicz inequalities of the type $\vert F(z)\vert\geq C\vert z\vert^\theta$ for families of polynomials defining a zero dimensional subspace of $\C^{d+1}$ were studied very early by P\l oski (see \cite{P1}). In the paper \cite{Hi2} Hickel uses methods of normalized blowing up and interpretation of $\overline\nu_{\Ii}(J)$ as the inverse of a \L ojasiewicz exponent to prove:
\begin{theorem}{\rm (Hickel)} Let $k$ be an algebraically closed field with a non trivial absolute value $\vert\ \vert$. Let $I=(p_1,\ldots ,p_m)\subset k[X_1,\ldots,X_n]$ and $V=\{z\in k^n/p_1(z)=\ldots = p_m(z)=0\}$. Set $d_i={\rm deg}p_i$. There exist a constant $C>0$ and exponents $\theta_1,\theta_2$ such that for all $z\in k^n$ the following inequality holds:
$$\sum_{i=1}^m\frac{\vert p_i(z)\vert}{(1+\vert z\vert)^{d_i}}\geq C\frac{{\rm d}(z,V)^{\theta_1}}{(1+\vert z\vert)^{\theta_2}},$$
where the best exponents $\theta_1,\theta_2$ are rational numbers, which can be computed from invariants which appear in the normalized blowing up of the sheaf of ideals on $\P^n(k)$ corresponding to the homogeneizations of the polynomials $p_i$.\end{theorem}
Here ${\rm d}(z,V)$ is the distance from $z$ to $V$. Of course Hickel's statement is much more precise, and the techniques of proof use more global techniques than just local intersection theory, in particular the refined B\'ezout theorem.  The close connexion between such inequalities and the effective Nullstellensatz has been known since the paper \cite{B} of Brownawell, and the paper of Hickel makes important headway towards optimal bounds for $s$ and the degrees of the $q_kf_k$ in terms of $n$ and the degrees of the polynomials $f_1,\ldots ,f_t$ and $g\in\sqrt{(f_1,\ldots ,f_t)k[x_1,\ldots,x_n]}$ in equalities such as $g^s=\sum_kq_kf_k$,  in all characteristics. See also \cite{L-T}, Compl\'ement 5, and the article \cite{C-K-T}.
\end{subsection}
\begin{subsection}{\L ojasiewicz exponent and log canonical threshold}
The gradient\break \L ojasiewicz exponents of the restrictions to plane sections of all dimensions appear in an estimate of the \emph{log canonical threshold} associated to the hypersurface singularity $f(z_1,\ldots, z_{d+1})=0$ at the origin. I refer to the notes \cite{Mu} of Mus\c tat\u a and the book \cite{La}, Vol. II, of Lazarsfeld for the definition. It is of the nature of a \L ojasiewicz exponent for differential forms.\par
The jacobian Newton polygon of a hypersurface with isolated singularity, say at the origin, can be written in the form:
$$\Newton_{j(f)}(m)=\sum_{q=1}^t\Teiss{e_q}{m_q},$$
where the $q$'s index the branches of a general polar curve, which means the critical locus of a map $F\colon \C^{d+1}\to\C^2,\ z\mapsto (\ell(z),f(z))$ with $\ell$ a general linear form. The integers $e_q,m_q$ are then computed from a parametrization of each branch: $m_q$ is the multiplicity of that branch, and $e_q+m_q$ is its intersection multiplicity with $f=0$. See \cite{T2}, where it is proved that the best possible exponents for the gradient \L ojasiewicz inequalities to hold in a neighborhood of the singular point,
$$\vert f(z)\vert^{\theta_1}\leq C_1\vert{\rm grad}f(z)\vert \ \ {\rm and} \ \ \vert z\vert^{\theta_2}\leq C_2\vert{\rm grad}f(z)\vert ,$$
are $\theta_1=\frac{\eta}{1+\eta}$ and $\theta_2=\eta$ where $\eta={\rm sup}_q\frac{e_q}{m_q}$. This shows not only that the \L ojasiewicz exponents are rational numbers and that $\theta_1<1$, a well known and very useful result of \L ojasiewicz in the real case (see \cite{L-T}, Compl\'ement 1), but also that they are invariants of Whitney equisingularity, since the jacobian Newton polygon is. \par
Since the jacobian Newton polygon of a general hyperplane section is well defined, we may for each $1\leq i\leq d+1$ define $\eta^{(i)}$ to be the gradient \L ojasiewicz exponent  and $\tau^{(i)}$ to be the ${\rm inf}_q\frac{e^{(i)}_q}{m^{(i)}_q}$ for the restriction of $f$ to a general $i$-dimensional linear space through the origin. They are respectively the  largest  and smallest inclinations of the edges of the jacobian Newton polygon of that restriction. We agree that the restriction of $f$ to a general line, which is equivalent to $z_1^m$, gives $\theta^{(1)}=\tau^{(1)}=m-1$ where $m$ is the order of $f$ at $x$, the multiplicity at $x$ of the hypersurface $f=0$. Then we have:
\begin{theorem}{\rm (Loeser, see \cite{Lo1}, \cite{Lo2})} Assume that the isolated singularity of $f(z)=0$ at $x$ is not a rational singularity (this is equivalent to ${\rm lct}_x(f)\leq 1$ by a result of Morihiko Saito). The log canonical threshold ${\rm lct}_x(f)$ of $f$ at $x$ satisfies the inequalities
$$\sum_{i=1}^{d+1}\frac{1}{1+\lceil \eta^{(i)}\rceil}\leq {\rm lct}_x(f)\leq \sum_{i=1}^{d+1}\frac{1}{1+\lfloor \tau^{(i)}\rfloor}\ .$$
\end{theorem}
In the case where $f=z_1^{a_1}+\cdots+z_{d+1}^{a_{d+1}}$ with $a_1\leq a_2\leq\cdots \leq a_{d+1}$, one computes (see \cite{T2}) that $\tau^{(i)}=\eta^{(i)}=a_i-1$ and recovers the well known formula ${\rm lct}_0(f)=\sum_{i=1}^{d+1}\frac{1}{a_i}$ when this sum is $\leq 1$. \par\noindent
In \cite{Lo2}, one finds an application to arithmetical algebraic geometry. I had in \cite{T3} conjectured the first inequality, without the $\lceil\ \rceil$, for reasons related to the rate of vanishing as $\lambda$ tends to zero of the volumes of the vanishing cycles on a Milnor fiber $(f(z)=\lambda)\cap\B(x,\epsilon)$.\par
It is shown in \cite{T2} that:\par\noindent
$\bullet$ The length of the projection to the horizontal (resp. vertical) axis of the jacobian Newton polygon $\Newton_{j(f\vert H^{(i)})}(m)$ of the restriction of $f$ to a general vector subspace of dimension $i+1$ is equal to the Milnor number $\mu^{(i+1)}(f)$ (resp. $\mu^{(i)}(f)$) of the restriction of $f$ to $H^{(i+1)}$ (resp $H^{(i)}$).  This jacobian Newton polygon is therefore a refinement of the elementary polygon $\Teiss{\mu^{(i+1)}(f)}{\mu^{(i)}(f)}$.\par\noindent
$\bullet$
The double of the \emph{area} bounded by the jacobian Newton polygon is the Milnor number of a general hyperplane section of the Thom-Sebastiani "double" $f(z_1,\ldots, z_{d+1})+f(w_1,\ldots, w_{d+1})$ of the function $f(z_1,\ldots, z_{d+1})$.\par\noindent

\begin{remark} If $f$ is a homogeneous polynomial of degree $m$ with isolated singularity, one finds that all $\frac{e_q}{m_q}$ are equal to $m-1$, so that the
optimal inequality is the expected ''homogeneous'' one:
$$\vert f(z)\vert^{m-1}\leq C\vert \hbox{\rm grad}f(z)\vert^m.$$
But this is no longer true if the cone defined by the vanishing of the homogeneous polynomial $f$ has a non isolated singularity at the origin; here is a
counterexample: $$f(z_1,z_2,z_3)=z_3(z_1z_3^3+z_2^4)^2+z_1^9.$$
The plane projective curve determined by this equation has at the point $z_1=z_2=0$ a singularity isomorphic to $w^2-z^{36}=0$; it corresponds to a number
$\eta=\hbox{\rm sup}_q\frac{e_q}{m_q}$ equal to $35$, which must manifest itself in any neighborhood of the vertex of the cone defined by the equation,
and prevents the exponent there  from being $\frac{8}{9}$, since it must be at least equal to $\frac{35}{36}$. \emph{(This is an example I produced in the late 1980's in answer to a question of Amnon Neeman; the question of finding a sharp bound, in this homogeneous case, for $\theta_1$ in terms of $m$ is still open as far as I know)}.
\end{remark}
\end{subsection}

\section{The gradient \L ojasiewicz exponent in dimension 2}
In \cite{GB-P} and \cite{GB-K-P} some unexpected results, which find their roots in an earlier work of A. P\l oski, were proved, concerning the rational numbers which can appear as $\theta_2$, in the notations of the previous subsection, in the case $d=1$ i.e., for plane curve singularities. 
\begin{theorem}{\rm (P\l oski, see \cite{P2})} For $d=2$ the rational number $\theta_2$ belongs to the set $$\{N+\frac{b}{a}\ , N,a,b\in \Z,\ 0\leq b<a<N\}.$$
\end{theorem}
\begin{theorem}{\rm (Garc\'ia Barroso-P\l oski, see \cite{GB-P})} The elements of the preceding set for which $a=N-1, b>1, {\rm gcd}(a,b)=1$ are \emph{not} exponents $\theta_2$ for a plane curve. 
\end{theorem}
There are more results of this kind in \cite{GB-K-P}.\par\noindent
\begin{remark}
The problem of determining from intrinsic invariants whether a given ideal is the jacobian ideal of a hypersurface is part of the folklore, and nobody knew how to begin. The previous theorem provides a partial answer in dimension $2$: if a primary ideal $I$ in $\C\{z_1,z_2\}$ is such that $(\overline\nu_I(m))^{-1}$ belongs to the set determined in  \cite{GB-P}, then it cannot be a jacobian ideal.
\end{remark}
A theorem of M. Merle in \cite{Me} implies that the jacobian Newton polygon of a plane branch is a total invariant of its equisingularity type: it determines and is determined by the Puiseux exponents, or the semigroup of values, or the topological type of its embedding in $\C^2$ locally at $x$. This is no longer true for   reduced plane curves, as was shown by Eggers (see \cite{E}), and the correct generalization was achieved by E. Garc\'ia Barroso in \cite{GB}.\par\medskip
If one thinks that the jacobian Newton polygon is a natural extension of the gradient \L ojasiewicz exponent it is natural to ask whether all Newton polygons meeting both axes (or \emph{convenient}) can be the jacobian Newton polygon of a hypersurface with isolated singularity, for a given dimension $d$.\par\noindent
For $d=1$ there is a very nice answer, for branches, due to E. Garc\'ia Barroso and J. Gwo\'zdziewicz in \cite{GB-G}: there is an arithmetic characterization of the Newton polygons that are the jacobian Newton polygon of a plane branch. Moreover, they prove that if the jacobian Newton polygon of a plane curve is that of a branch, then the curve is a branch, giving rise to new criteria of irreducibility.
\section{\L ojasiewicz exponents in infinite dimensions}
Let us briefly go back to a \L ojasiewicz inequality in the real case, relating the distance to the zero set of a function to the value of the function:
\begin{theorem} \emph{(\L ojasiewicz, \cite{L})} Let $\Omega$  be an open subset of $\R^n$, let $f$ be an
analytic function on $\Omega$ and $V$ the locus of zeroes of $ f$ . For every compact
$K\subset \Omega$ there exist $ C>0, \alpha \geq 0$ such that:
$$\forall x\in K,\ \hbox{\rm dist}(x,V)\leq C\vert f(x)\vert^\alpha ,$$\noindent
where $d(x,V)$ is the distance from $x$ to $V$ in $\R^n$. If $V=\emptyset$,
we take $\alpha =0$ and agree that $d(x, \emptyset )^0=1$.\end{theorem}
The key to understanding why this \L ojasiewicz inequality is true is that the distance function to an analytic set behaves much like an analytic function: it is subanalytic. Note that in comparison with the previous inequalities, we have $\alpha=\theta^{-1}$.\par\noindent
Let us reformulate the statement of the theorem as follows: for each
sufficiently small $\epsilon$,\par\noindent
$$\hbox{\rm if }\ \vert f(x)\vert \leq\epsilon, \hbox{\rm there exists }x_0\in V \ \hbox{\rm such that }\ f(x_0)=0\ \hbox{\rm and }\vert x-x_0\vert\leq
C\epsilon^{\alpha}.$$
Let us now make a translation which is perhaps "classical" (see the review of \cite{Sc}) and certainly well known to Hickel (see \cite{Hi1}, 2.2): consider a space of analytic or formal functions in variables $x_1,\ldots ,x_n$, for example the maximal ideal ${\mathcal M}$ of the ring ${\O}={\bf
C}\{x_1,\ldots ,x_n\}$ of convergent power series at the origin. Consider a finite set of algebraic (or analytic) equations $\Phi_k(x_1,\ldots ,x_n,Y_1,\ldots
,Y_p)=0$, with $\Phi_k\in {\mathcal O}[Y_1,\ldots ,Y_p]$ (or in ${\mathcal O}\{Y_1,\ldots ,Y_p\}$) and $0\leq k\leq q$, which we shall forget to write from now
on; they define a subset
$V$ of the infinite dimensional space
${\mathcal M}^p$, which we shall call an {\it algebraic subset} (resp. an {\it analytic subset}) of this infinite-dimensional space.\par\noindent
On ${\mathcal M}$ there is a metric given by the norm $\Vert g\Vert=e^{-\nu_{\mathcal M}(g)}$, where $\nu_{\mathcal M}(g)$ is the ${\mathcal M}$-adic order. We have
$$\Vert g+g'\Vert \leq \hbox{\rm max}(\Vert g\Vert ,\Vert g'\Vert).$$
We extend this norm to ${\mathcal M}^p$ by taking the maximum of the norms of the coordinates.\par\noindent
In this context, \L ojasiewicz's inequality becomes, after a small translation as above: \par\noindent
For any sufficiently large integer $N$, if we have series without constant term ${\mathbf y}^{(N)}(x)=( y^{(N)}_1(x),\ldots , y^{(N)}_p(x))\in {\mathcal M}^p$ such that for all $k$ the
inequalities $\nu_{\mathcal M}(\Phi_k(x,{\mathbf y}^{(N)}(x)))\geq N$ hold, there exists a point ${\mathbf y}(x)\in V$ with $\Vert{\mathbf y}^{(N)}(x)-{\mathbf y}(x) \Vert\leq
Ce^{-N\alpha}$. If we write the positive constant $C$ in the form $C=e^{-b}$ this
becomes, with $a=\alpha\in {\bf R}_+,\ b\in {\bf R}$:
$$\nu_{\mathcal M}({\mathbf y}^{(N)}(x)-{\mathbf y}(x))\geq aN+b,$$
which means that each coordinate is in ${\mathcal M}^{aN+b}$, so that we can reinterpret the \L ojasiewicz inequality:
\begin{proposition}\label{Lojinf} Let $V\subset {\mathcal M}^p$ be an algebraic subset
defined by the polynomials $\Phi_k\in {\mathcal O}[Y_1,\ldots ,Y_p],\ 1\leq k\leq q$. The following statements are equivalent:\par\noindent
a) The subset $V$ and its defining equations satisfy \L ojasiewicz's inequality with respect to the norm $\Vert f\Vert$.\par\noindent 
b) There exist  constants
$a\in {\mathbf R}_+,\ b\in {\mathbf R}$ depending only on the system of equations
$\Phi_k$ such that for each sufficiently large integer
$N$, if there exits a system of series without constant terms ${\mathbf y}^{(N)}(x)\in {\mathcal M}^p$ such that for all $k$ we have $\nu_{\mathcal M}(\Phi_k(x, {\mathbf y}^{(N)}(x)))\geq N$, then there exists a system
${\mathbf y}(x)\in {\mathcal M}^p $ such that $\Phi_k(x,{\mathbf y}(x))=0$ and 
$\nu_{\mathcal M}({\mathbf y}^{(N)}(x)-{\mathbf y}(x))\geq aN+b$.\end{proposition}
It is natural to ask \emph{whether the \L ojasiewicz inequality of Proposition \ref{Lojinf} holds for
algebraic and analytic subsets of ${\mathcal M}^p$}.\par Strong supporting evidence comes from M. Artin's
theorem:
\begin{theorem} {\rm (M. Artin, \cite{A})} Given a system of equations such as $\Phi_k(x,y)$, there exists an application $\beta\colon {\bf N}\to{\bf N}$ such that for
any integer $i$, if there exists a system of series without constant term  $\tilde y(x)$ such that $$\nu_{\mathcal M}(\Phi_k(x,\tilde y(x)))\geq \beta (i),$$ then there
exists a system $y_0(x)$ such that $\Phi_k(x,y_0(x))=0$ and $$\nu_{\mathcal M}(\tilde y(x)-y(x))\geq i.$$\end{theorem}
The only difference is that the function $\tilde\beta (i)=\frac{i-b}{a}$ obtained by inverting the map $N\mapsto aN+b$ of statement \ref{Lojinf} is linear, which increases the interest of finding linear upper bounds for Artin's beta function. 
The following result due to Hickel may help to throw some light and make one optimistic about the linearity:\par
Let our indeterminates consist of one indeterminate $t$ and let the equation $\Phi$ be $f(Y_1,\ldots ,Y_n)=0$, where $f(z_1,\ldots ,z_n)$ is an analytic function
with an isolated critical point at the origin. Solving the equation $f(y_1(t),\ldots ,y_n(t))=0$ with series without constant terms amounts to finding analytic arcs
through the origin on the hypersurface defined by $f(z_1,\ldots ,z_n)=0$. Keeping the notation $\eta=\hbox{\rm sup}_q\frac{e_q}{m_q}$ introduced above,
we have:\par\noindent
\begin{theorem}{\rm (Hickel, \cite{Hi1})} The function $\beta$ of Artin satisfies in this case $\beta(i)\leq \lfloor\eta i\rfloor+i$.\end{theorem}
However, the existence in general of a linear upper bound for Artin's beta function, conjectured by Spivakovsky in [S], was shown to be false by G. Rond in \cite{Ro1}.\par\noindent
So the nature of \L ojasiewicz inequalities changes when one passes to infinite dimensional spaces, if we insist on taking the metric $\Vert g\Vert=e^{-\nu_{\mathcal M}(g)}$.
\par\medskip\noindent
A recent result of Moret-Bailly shows that if we seek roots of our equations\break $\Phi_k(x_1,\ldots ,x_n,Y_1,\ldots
,Y_p)=0$ not only in ${\mathcal M}$ but in the maximal ideal ${\mathcal M}_\nu$ of a henselian valuation ring $R_\nu$ containing $\O$ and such that ${\mathcal M}_\nu\cap\O={\mathcal M}$, and use an adapted metric, we recover a \L ojasiewicz inequality. The following statement is the special case for an affine $R_\nu$-scheme, of Moret-Bailly's theorem:\par\noindent
Let $R_\nu$ be a henselian valuation ring such that the field of fractions $\hat K$ of its completion is a separable extension of its field of fractions $K$. Let $\Gamma$ be the group of values of the valuation $\nu$, which is totally ordered and, for $\gamma\in \Gamma_+$, consider the valuation ideal $\Pp_\gamma(R_\nu)=\{x\in R_\nu\vert \nu(x)\geq \gamma\}$.

\begin{theorem}{{\rm (Moret-Bailly, \cite{MB})}} Consider a finite set of algebraic equations\break  $\Phi_k(Y_1,\ldots
,Y_p)=0$, with $\Phi_k\in R_\nu[Y_1,\ldots ,Y_p]$. There exist a positive integer $B$ and an element $\delta\in \Gamma_+$ with the following properties: for each $\gamma\in \Gamma_+$ and ${\mathbf y}^{(\gamma)}=(y_1^{(\gamma)},\ldots, y_p^{(\gamma)})\in R_\nu^p$ such that $\Phi_k(y_1^{(\gamma)},\ldots
,y_p^{(\gamma)})\in \Pp_{B\gamma+\delta}(R_\nu)$ for all $k$, there exists ${\mathbf y}\in R_\nu^p$ such that $\Phi_k({\mathbf y})=0$ for all $k$ and ${\mathbf y}^{(\gamma)}-{\mathbf y}\in (\Pp_\gamma (R_\nu))^p$.
\end{theorem}
So, in this situation, the Artin function associated to the valuation $\nu$ (and not to $\nu_{\mathcal M}$) is indeed bounded by the linear function $\gamma\mapsto B\gamma+\delta$.\par\noindent
If we restrict ourselves to valuations of rank one, for which the value group is an ordered subgroup of $\R$, and take the metric $\Vert x\Vert=e^{-\nu(x)}$ we are exactly in the same situation as at the beginning of this section. Indeed, this suggests that we must seek solutions in henselian valuation rings, and that the right "distance to the zero set" is measured by the valuation. Moret-Bailly's theorem gives us a version of the \L ojasiewicz distance inequality in the space $R_\nu^p$.

%-----------------------------------------------

\par\vskip1truecm
\centerline{copyright Bernard Teissier, 2012}
\end{document}